\documentclass[fleqn,11pt]{article}
\usepackage{amssymb}
\usepackage{graphicx}
\usepackage{hyperref}
\usepackage{amsmath}
\usepackage{times}
\begin{document}
\title{Fast multipole accelerated boundary element methods for room acoustics}
\author{Nail A. Gumerov\thanks{\tt {ngumerov@umd.edu}, {https://www.umiacs.umd.edu/users/gumerov}} \ and Ramani Duraiswami\thanks{\tt {ramanid@umd.edu}, {https://www.umiacs.umd.edu/users/ramani}} \\
	Perceptual Interfaces and Reality Lab, UMIACS,\\
	University of Maryland, College Park, MD\\  Also at VisiSonics Corporation, College Park,	MD}

\maketitle

\abstract{The direct and indirect boundary element methods, accelerated via the fast multipole method, are applied to numerical simulation of room acoustics for large rooms of volume $\sim 150$ $m^{3}$ and frequencies up to 5 kHz on a workstation. As the parameter $kD$ (wavenumber times room diameter) is large, stabilization of the previously developed FMM algorithms is required for accuracy. A stabilization scheme is one of the key contribution of this paper. The computations are validated using well-known image source solutions for shoebox shaped rooms. Computations for L-shaped rooms are performed to illustrate the ability to capture diffractions. The ability to model in-room baffles, and boundary openings (doors/windows) is also demonstrated. The largest case has $kD>1100$ with a discretization of size 6 million elements. The performance of different boundary integral formulations was compared, and their rates of convergence using a preconditioned flexible GMRES were found to be substantially different. These promising results suggest a path to efficient simulations of room acoustics via high performance boundary element methods.}

\section{Introduction}

The perception of sound by humans, and the computational simulation of
acoustic propagation leading to such perception, have been the subject of
many studies. The simulation models the propagation of sound from the source
location to the receiver -- a microphone or the ears of a listener. The
setting may be a real source in a real environment -- indoors or outdoors,
or a virtual environment -- in a game or extended reality simulation. The
simulation is often performed in two independent pieces -- first the
multipath propagation of the sound from the source to the location of the
receiver, and second the scattering of the incoming sound off the anatomy en
route to the ear canal of the listener.

The equations governing acoustic propagation is the linear wave equation
neglecting dispersion. Time invariant (or slowly varying) systems can be
usefully characterized by an impulse response, which can then be used with
convolution to compute the sound field for arbitrary sources. The scattering
off the anatomy can be characterized by the head related impulse response
(HRIR), and the scattering off the room walls and other obstacles by the
point-to-point room impulse response (RIR). The overall binaural room
impulse response (BRIR) can be obtained via a compositional model from these %
\cite{Zotkin2004:IEEE}.

This work is concerned with the computation of the room impulse response
using fast multipole accelerated boundary element (FMMBEM) methods. A
previous paper used FMMBEM to compute the HRIR \cite{Gumerov2010:JASA}. Once
the response is computed it can be put to various uses -- e.g., used in an 
\emph{audio engine} for rendering 3D audio, e.g., in a game or virtual
reality application \cite{Raghuvanshi2018:Siggraph}; or used to assess the
listening properties of a space such as a concert hall or a meeting room; or
assess, via simulation, architectural choices of room shapes and materials %
\cite{Kuttruff2017:book}.

Among methods for the solution of the wave equation are approximate
geometric methods that rely on tracing of rays, or groups of rays. Such
approaches are more popular (and more accurate for) light propagation but
can be adapted to acoustics \cite{Funkhouser2004:JASA, Savioja2015:JASA}.
These methods can be augmented with approaches that approximately account
for wave effects such as diffraction and diffuse reflections \cite%
{Schissler2014:Siggraph, Potter2020:ASA}. The methods can quickly get
leading order effects in an explainable fashion. However, when used to
obtain the whole impulse response, they suffer from combinatorial explosion
when the propagation process has many interactions of waves with surfaces,
and with each other (caustics). Moreover, it is not possible to provide an
accuracy bound or a computational cost bound for achieving a given numerical
accuracy. Further, it may not be possible to enforce realistic boundary
conditions. On the other hand, these methods allow for simple geometries to
be computed efficiently (e.g., via images \cite{Allen1979:JASA,
Borish1984:JASA}). Moreover, the methods allow taking large time steps as
the wavefront is propagated. Finally, because of the development of GPU
based raytracing approaches, there is promise for fast execution on modern
computers.

In contrast the methods that solve the wave equation can usually provide an
error and computational cost bounds, as well as the ability to enforce
correct boundary conditions. The wave equation can be solved in time and
space variables over the entire domain. This is done in discretization of
its original form via finite differences \cite{Hamilton2017:IEEE},
finite-volume\cite{Bilbao2016:IEEE} or spectral-element methods \cite%
{Raghuvanshi2009:TVCG}, or finite element methods \cite{Okuzono2016:AA}. A
key parameter that enters into the accuracy of the simulation is the
smallest wavelength that must be resolved (or equivalently the highest
frequency). The Nyquist criterion requires at least two discretization
points (or elements or degrees of freedom) per wavelength ($\lambda$), while
in practice this number is six to ten for accuracy. In three dimensions,
this means that if the average linear dimension is $D$, and the wavenumber $%
k $ ($=2\pi/\lambda$), the number of variables in the discretization is $%
N\simeq O(kD)^3$.

The goal of these simulations is to obtain the response due to an impulsive
source at a particular location. The duration of the simulation $Y$ is the
length of the impulse response -- the time taken for a unit impulse to die
down to negligible levels. When the equation is solved via a time marching
method, with an average time step size $h$, the number of steps is $T/h$.
When a linear system is solved at each step, the cost of the solution will
be proportional to $O(N)^\alpha$. For a direct solve via LU decomposition
the value of $\alpha$ is 3. If sparsity or other properties can be used to
control the cost of the solution, either via an iterative Krylov method or a
fast direct method, the value of $\alpha$ can be reduced to 2, or even lower
to 1 using specialized direct methods, which are an item of research. So,
the overall cost of computation can be reduced to $O(T/h)(kD)^{3\alpha}$.

These global methods have several other pros and cons. They can solve more
generalized versions of the wave equation (e.g., with dispersion, and
variable sound speeds). Moreover, once the solution is completed, the
pressure and velocity are known throughout the field. On the other hand, the
methods cannot handle infinite domains (with radiation conditions) easily;
though this is not necessary for internal acoustics problems. The treatment
of curved boundaries, or those with corners can lead to errors. Finally, it
is not straightforward to incorporate frequency dependent parameters in
boundary conditions.

The wave equation is often solved in the Fourier domain, which leads to the
Helmholtz equation for each frequency component, and the frequency response
is obtained instead of the impulse response. While, in principle, we need to
solve this equation over an infinite band of frequencies, we only need to
obtain the complex response at $T/2h$ frequencies to compute the impulse
response. While the Helmholtz equation can also be solved via domain
discretization using the finite difference or finite element methods, in
this paper we seek to solve it via the boundary element method. For the BEM
we only need to discretize the boundary of the domain leading to $N\simeq
O(kD)^{2}$ sized discretization, with dense matrices.

Note that for room acoustics the value of $kD$ can be relatively large --
from a few hundreds to tens of thousands, depending on the dimensions of the
room and the maximum frequency considered. Because the discretization size
varies as a power of this parameter (as the cube for domain methods, or as
the square for boundary methods), the resulting linear systems are very
large, and the time taken for setting up the linear system, as well as for
their solution is significant, unless fast algorithms are used. We explore
two variants of the boundary-element method, the direct \cite%
{Gumerov2009:JASA} and the indirect boundary element method \cite%
{Gumerov2013:POMA}, in Sections \ref{BEM} and \ref{IBEM}.

Using the fast-multipole method for accelerating matrix-vector products that
occur in Krylov iterative methods for solving the linear systems can
significantly reduce the time and memory needed for solution. For the high
values of $kD$ that are encountered in room acoustics, the translation and
rotation operators in the FMM \cite{Gumerov2009:JASA, Gumerov2010:JASA} were
found to be inaccurate due to instabilities in the recursion procedures that
are used. We propose in Section \ref{FMM} new stable translation and
rotation procedures that are stable at higher $kD$. Also, we compare the
integral equation formulations for their convergence rates, and find that
some formulations are more optimal. We achieve solutions with $O(kD)^2$
memory and $O(kD)^{2\alpha}$ time, with $1\leq \alpha < 2$. Results on
sample problems are shown, and possible directions for future work
indicated. Based on these results we believe FMM accelerated BEM is a
promising approach for accurate and fast room acoustics simulations for
practical room sizes and frequency ranges.

\section{Statement of the problem}

We consider a room bound by closed or open surface $S$ (the walls, ceiling,
and floor, possibly with some openings and some internal baffled surfaces).
Assume that the inhomogeneity of the air due to the temperature gradient is
negligible and the propagation of sound can be described by the wave
equation with a constant speed of sound, $c$,%
\begin{equation}
\frac{\partial ^{2}p\left( \mathbf{r},t\right) }{\partial t^{2}}=c^{2}\nabla
^{2}p\left( \mathbf{r},t\right) ,  \label{sp1}
\end{equation}%
where $p$ is the acoustic pressure depending on time $t$ and the
radius-vector $\mathbf{r}$ or a point in the room. Furthermore, we consider
time-periodic signals with period $T$. The Fourier transform in time domain
results in the Helmholtz equation for each harmonic component, 
\begin{equation}
\nabla ^{2}\psi \left( \mathbf{r};k\right) +k^{2}\psi \left( \mathbf{r}%
;k\right) =0,\quad k=\omega /c,  \label{sp2}
\end{equation}%
where $\psi $ is the complex amplitude of the sound pressure, $\omega $ is
the circular frequency. Assume further that the room boundaries are locally
reacting surface, whose properties do not depend on the angle of sound
incidence, so the following boundary conditions can be imposed \cite%
{Kuttruff2017:book} 
\begin{equation}
\left. -\zeta (\mathbf{r};k)\frac{\partial \psi }{\partial n}+ik\psi \right|
_{S}=0.\quad  \label{sp3}
\end{equation}%
The normal to $S$ is assumed directed outside the room, $\zeta $ is the
specific impedance, which is complex, and related to the wall reflection ($R$%
) and absorption ($\alpha )$ coefficients as%
\begin{equation}
\zeta (\mathbf{r};k)=\frac{1+R(\mathbf{r};k)}{1-R(\mathbf{r};k)},\quad
\alpha (\mathbf{r};k)=1-\left| R(\mathbf{r};k)\right| ^{2}.  \label{sp4}
\end{equation}%
For fully absorbing surfaces $\left( R=0\text{, }\alpha =1\right) $, $\zeta
=1$. For sound hard surfaces $\left( R=1\text{, }\alpha =0\right) $, $\zeta
=\infty $, boundary condition (\ref{sp3}) is replaced by the Neumann
condition%
\begin{equation}
\left. \frac{\partial \psi }{\partial n}\right| _{S}=0.  \label{sp5}
\end{equation}
Due to linearity of the problem it is sufficient to obtain a solution for a
single monopole source of unit intensity located at some point $\mathbf{r}%
_{s}$ as more complex incident fields can be simulated by integration over
prescribed source distributions. Thus, we consider 
\begin{equation}
\psi (\mathbf{r};k)=\psi ^{in}(\mathbf{r};k)+\psi ^{sc}(\mathbf{r};k),\quad
\psi ^{in}=G_{k}\left( \mathbf{r},\mathbf{r}_{s}\right) =\frac{e^{ik\left| 
\mathbf{r}-\mathbf{r}_{s}\right| }}{4\pi \left| \mathbf{r}-\mathbf{r}%
_{s}\right| },  \label{sp6}
\end{equation}%
where $\psi ^{in}$ and $\psi ^{sc}$ are the incident and the scattered
fields, and $G_{k}\left( \mathbf{r},\mathbf{r}_{s}\right) $ is the Green's
function for the Helmholtz equation with wavenumber $k$. The field $\psi
^{sc}$ has no singularity inside the room and can be computed using boundary
conditions%
\begin{equation}
\left. \frac{\partial \psi ^{sc}}{\partial n}+\frac{k}{i\zeta (\mathbf{r};k)}%
\psi ^{sc}\right| _{S}=\gamma (\mathbf{r};k)=-\left. \left( \frac{\partial
\psi ^{in}}{\partial n}+\frac{k}{i\zeta (\mathbf{r};k)}\psi ^{in}\right)
\right| _{S},  \label{sp7}
\end{equation}%
Note that a solution $\psi$ obtained for such a formulation satisfies in the
domain%
\begin{equation}
\nabla ^{2}\psi \left( \mathbf{r};k\right) +k^{2}\psi \left( \mathbf{r}%
;k\right) =-\delta \left( \mathbf{r},\mathbf{r}_{s}\right).  \label{sp8}
\end{equation}

\section{Method of solution}

\subsection{Direct BEM}

\label{BEM} When the room surface $S$ is closed, the direct boundary element
method based on the use of Green's identity for internal points can be used, 
\begin{equation}
\psi ^{sc}(\mathbf{r,}k)=L_{k}\left[ \frac{\partial \psi ^{sc}}{\partial n}%
\right] \left( \mathbf{r}\right) -M_{k}\left[ \psi ^{sc}\right] \left( 
\mathbf{r}\right) ,\quad \mathbf{r}\notin S,  \label{dbem1}
\end{equation}%
where $L_{k}$ and $M_{k}$ are the single and double layer potentials, 
\begin{eqnarray}
L_{k}\left[ \sigma \right] \left( \mathbf{r}\right) &=&\int_{S}\sigma \left( 
\mathbf{r}^{\prime }\right) G_{k}\left( \mathbf{r},\mathbf{r}%
^{\prime}\right) dS\left( \mathbf{r}^{\prime }\right) ,  \label{dbem2} \\
M_{k}\left[ \mu \right] \left( \mathbf{r}\right) &=&\int_{S}\mu \left( 
\mathbf{r}^{\prime }\right) \frac{\partial }{\partial n\left( \mathbf{r}%
^{\prime }\right) }G_{k}\left( \mathbf{r},\mathbf{r}^{\prime }\right)
dS\left( \mathbf{r}^{\prime }\right) .  \notag
\end{eqnarray}%
The surface values of the scattered field can be found by applying the
limiting form of Green's identity (\ref{dbem1}), which for smooth surfaces
yields%
\begin{equation}
\frac{1}{2}\psi ^{sc}(\mathbf{r,}k)=L_{k}\left[ \frac{\partial \psi ^{sc}}{%
\partial n}\right] \left( \mathbf{r}\right) -M_{k}\left[ \psi ^{sc}\right]
\left( \mathbf{r}\right) ,\quad \mathbf{r}\in S.  \label{dbem3}
\end{equation}%
Here the integral $M_{k}$ is singular, and should be treated in the sense of
principal value.

For room acoustics the boundary integral equation (BIE) (\ref{dbem3})
combined with the boundary condition (\ref{sp7}) is sufficient to solve the
problem. There is no need to introduce the treatment to avoid spurious
resonances, that was proposed by \cite{Burton1971:PRS} for the external
problem. The degeneracy of the BIE in that case is due to two reasons:
first, that the surface is sound hard, and, second that the problem is
external. For impedance boundary conditions with non-zero real part of $%
\zeta $ (\ref{sp4}), these resonances fo not arise. Even if the surface is
sound-hard, the problems in the BIE solution is related to the resonance in
the internal problem and is a numerical artifact. However, for internal
domains such resonances are physical and should be handled by different
means.

For impedance boundary conditions (\ref{sp7}) the BIE can be written either
for unknown $\psi ^{sc}$,%
\begin{equation}
\frac{1}{2}\psi ^{sc}(\mathbf{r,}k)+M_{k}\left[ \psi ^{sc}\right] \left( 
\mathbf{r}\right) +kL_{k}\left[ \frac{1}{i\zeta }\psi ^{sc}\right] \left( 
\mathbf{r}\right) =L_{k}\left[ \gamma \right] \left( \mathbf{r}\right)
,\quad \mathbf{r}\in S,  \label{dbem4}
\end{equation}%
or unknown $\partial \psi ^{sc}/\partial n$,%
\begin{equation}
\frac{1}{2}\frac{\partial \psi ^{sc}}{\partial n}(\mathbf{r,}k)+\frac{1}{%
\zeta (\mathbf{r},k)}M_{k}\left[ \zeta \frac{\partial \psi ^{sc}}{\partial n}%
\right] \left( \mathbf{r}\right) +\frac{k}{i\zeta (\mathbf{r},k)}L_{k}\left[ 
\frac{\partial \psi ^{sc}}{\partial n}\right] \left( \mathbf{r}\right) =%
\frac{1}{2}\gamma \left( \mathbf{r,}k\right) +\frac{1}{\zeta (\mathbf{r},k)}%
M_{k}\left[ \zeta \gamma \right] \left( \mathbf{r}\right) ,\quad \mathbf{r}%
\in S.  \label{dbem5}
\end{equation}%
If $\zeta $ is constant over frequency, the matrices, $A_{k}=\frac{1}{2}%
I+M_{k}$ $+\frac{k}{i\zeta }L_{k}$, are the same in both formulations, while
the right hand sides, $b=L_{k}\left[ \gamma \right] $ and $b=\left( \frac{1}{%
2}I+M_{k}\right) \left[ \gamma \right] $ are different (here $I$ is the
identity matrix). In our numerical experiments we observed faster
convergence for the flexible GMRES procedure for Eq. (\ref{dbem4}) and slow
or no convergence for Eq. (\ref{dbem5}) for large values of $\zeta $, while
the results, for convegrent cases, were correct for both formulations. The
direct BEM used in this paper is described in \cite{Gumerov2009:JASA}.

\subsection{Indirect BEM}

\label{IBEM} The indirect BEM can be used wherever the direct BEM can, and
in addition for cases when the surface $S$ is open (say a problem with an
open door to the room, and with some thin baffles inside the room). The
reason, why the direct BEM is preferable for closed rooms is that it leads
to a more compact BIE. For open surfaces we must use a representation of the
sound field via single and double layer potentials of unknown intensities $%
\sigma $ and $\mu $%
\begin{equation}
\psi ^{sc}(\mathbf{r};k)=L_{k}\left[ \sigma \right] \left( \mathbf{r}\right)
+M_{k}\left[ \mu \right] \left( \mathbf{r}\right) ,\quad \mathbf{r}\notin S.
\label{ibem1}
\end{equation}%
These intensities can be determined from the boundary conditions and the
jump conditions on the surface. Denote the side of $S$ facing the inside as $%
S^{-}$ and the side faced outside the room as $S^{+}$. Using ``+'' and ``-''
superscripts for parameters on that surfaces we can write the jump
conditions for field (\ref{ibem2}) as%
\begin{eqnarray}
\psi ^{sc\pm }(\mathbf{r};k) &=&L_{k}\left[ \sigma \right] \left( \mathbf{r}%
\right) +M_{k}\left[ \mu \right] \left( \mathbf{r}\right) \pm \frac{1}{2}\mu
\left( \mathbf{r},k\right) ,\quad \mathbf{r}\in S^{\pm },  \label{ibem2} \\
\frac{\partial }{\partial n}\psi ^{sc\pm }(\mathbf{r};k) &=&L_{k}^{\prime }%
\left[ \sigma \right] \left( \mathbf{r}\right) +M_{k}^{\prime }\left[ \mu %
\right] \left( \mathbf{r}\right) \mp \frac{1}{2}\sigma \left( \mathbf{r}%
,k\right) ,\quad \mathbf{r}\in S^{\pm },  \notag \\
L_{k}^{\prime }\left[ \sigma \right] \left( \mathbf{r}\right)
&=&\int_{S}\sigma \left( \mathbf{r}^{\prime }\right) \frac{\partial }{%
\partial n\left( \mathbf{r}\right) }G_{k}\left( \mathbf{r},\mathbf{r}%
^{\prime }\right) dS\left( \mathbf{r}^{\prime }\right) ,  \label{ibem3} \\
M_{k}^{\prime }\left[ \mu \right] \left( \mathbf{r}\right) &=&\frac{\partial 
}{\partial n\left( \mathbf{r}\right) }\int_{S}\mu \left( \mathbf{r}^{\prime
}\right) \frac{\partial }{\partial n\left( \mathbf{r}^{\prime }\right) }%
G_{k}\left( \mathbf{r},\mathbf{r}^{\prime }\right) dS\left( \mathbf{r}%
^{\prime }\right) ,  \notag
\end{eqnarray}%
and the integral singularities for evaluation points on the surface are
treated in terms of principal values. The normal is directed from ``-'' to
`+'' side of the surface.

Different variables can be declared as principal unknowns. For example, jump
relations (\ref{ibem2}) and boundary conditions (\ref{sp7}) allow us to
express the system in terms of two variables $\psi ^{sc+}$ and $\psi ^{sc-}$
on sides $S^{+}$ and $S^{-}$, respectively, as we have (we drop arguments $%
\left( \mathbf{r},k\right) $ for all functions to shorten notation)%
\begin{eqnarray}
\mu &=&\psi ^{sc+}-\psi ^{sc-},  \label{ibem4} \\
\sigma &=&\frac{\partial }{\partial n}\psi ^{sc-}-\frac{\partial }{\partial n%
}\psi ^{sc+}=\gamma ^{-}+\gamma ^{+}-\frac{k}{i\zeta ^{-}}\psi ^{sc-}-\frac{k%
}{i\zeta ^{+}}\psi ^{sc+},  \notag \\
\pm \frac{\partial }{\partial n}\psi ^{sc\mp }+\frac{k}{i\zeta ^{\mp }}\psi
^{sc\mp } &=&\gamma ^{\mp }=\mp \frac{\partial }{\partial n}\psi ^{in\mp }-%
\frac{k}{i\zeta ^{\mp }}\psi ^{in\mp },  \notag
\end{eqnarray}%
So, we have from Eq. (\ref{ibem2})%
\begin{gather}
\frac{\left( \psi ^{sc-}+\psi ^{sc+}\right)}{2} -M_{k}\left[ \psi
^{sc+}-\psi ^{sc-}\right] +L_{k}\left[ \frac{k\psi ^{sc-}}{i\zeta ^{-}}+%
\frac{k\psi ^{sc+}}{i\zeta ^{+}}\right] =L_{k}\left[ \gamma ^{-}+\gamma ^{+}%
\right] ,  \label{ibem5} \\
-M_{k}^{\prime }\left[ \psi ^{sc+}-\psi ^{sc-}\right] +\frac{1}{2}\left( 
\frac{k\psi ^{sc+}}{i\zeta ^{+}}-\frac{k\psi ^{sc-}}{i\zeta ^{-}}\right)
+L_{k}^{\prime }\left[ \frac{k}{i\zeta ^{-}}\psi ^{sc-}+\frac{k}{i\zeta ^{+}}%
\psi ^{sc+}\right]  \notag \\
=\frac{1}{2}\left( \gamma ^{+}-\gamma ^{-}\right) +L_{k}^{\prime }\left[
\gamma ^{-}+\gamma ^{+}\right] .  \notag
\end{gather}%
Note that in contrast to the BIE equations (\ref{dbem3}), the equations (\ref%
{ibem5}) contain integrals $L_{k}^{\prime }$ and the hypersingular integral $%
M_{k}^{\prime }$. In our conference communication \cite{Gumerov2013:POMA},
we briefly described the indirect BEM from more general case than room
acoustics. The details for computing the boundary integrals can be found in %
\cite{Gumerov2021:Arxiv}.

\subsection{Fast multipole method}

\label{FMM} The FMM strategy for boundary element methods was discussed in %
\cite{Gumerov2009:JASA}. According to that approach, the discretization
should be selected in a way that approximation for far field integrals is
good enough. The space partitioning with octree should be designed in a way
that the size of the smallest box $\delta $ is larger than the cutoff
distance, $\delta >r_{c} $. In this computation of matrix entries
corresponding to the interactions between closely located elements can be
done using analytical formulae, while far field interactions are obtained
via regular FMM procedures involving multipole and local expansions and
translations.

The FMM for the Helmholtz kernel used in the present study is the same as in %
\cite{Gumerov2009:JASA}. It requires local-to-local (which are the same as
multipole-to-multipole) translations for the Helmholtz equation even at
large frequencies. Such translations are performed using the
RCR-decomposition of translations (rotation, coaxial translation, and back
rotation). We found that the recurrences for computation of the reexpansion
coefficients proposed in our paper \cite{Gumerov2004:SISC} and book \cite%
{Gumerov2004:book} are not stable at large values of parameter $kD$, where $%
D $ is the size of the computational domain (diagonal of the bounding cube).
In our room computations presented in this paper $kD$ may the value of
several thousands (which is also the order of magnitude for the degree of
spherical harmonics). This relates both to rotations and coaxial
translations. A conditionally stable recursive procedure for computation of
the rotation operators, which is tested up to degrees of spherical harmonics
10$^{4}$ is presented in our study \cite{Gumerov2015:inBook}. There exist
other rotation methods tested for degrees of spherical harmonics of the
order of 10$^{3}$ \cite{Gimbutas2009:JCP}, which potentially can be used in
the present method. Stabilization of the recursions for the coaxial
translation coefficients is new.

\subsubsection{Stable high frequency coaxial translations}

In the FMM two types of representation of the solutions of the Helmholtz
equation are considered, the local and the multipole (far field expansions).
These expansions are 
\begin{eqnarray}
\psi \left( \mathbf{r}\right)
&=&\sum_{n=0}^{p-1}\sum_{m=-n}^{n}C_{n}^{m}S_{n}^{m}\left( \mathbf{r}\right)
,\quad \psi \left( \mathbf{r}\right)
=\sum_{n=0}^{p-1}\sum_{m=-n}^{n}D_{n}^{m}R_{n}^{m}\left( \mathbf{r}\right) ,
\label{shft0} \\
S_{n}^{m}\left( \mathbf{r}\right) &=&h_{n}\left( kr\right) Y_{n}^{m}\left( 
\mathbf{r/}r\right) ,\quad R_{n}^{m}\left( \mathbf{r}\right) =j_{n}\left(
kr\right) Y_{n}^{m}\left( \mathbf{r/}r\right) ,\quad r=\left| \mathbf{r}%
\right| ,  \notag
\end{eqnarray}%
where $S_{n}^{m}\left( \mathbf{r}\right) $ and $R_{n}^{m}\left( \mathbf{r}%
\right) $ are the singular and regular at $r=0$ spherical basis functions, $%
j_{n}$ and $h_{n}$ are the spherical Bessel and Hankel functions of the
first kind, and $Y_{n}^{m}\left( \mathbf{r/}r\right) $ are the spherical
harmonics (we use definition which can be found in references \cite%
{Gumerov2004:book, Gumerov2009:JASA, Gumerov2015:inBook}). The infinite
summation over $n$ is truncated to $p$, leading to $p^{2}$ coefficients.
Each step of the FMM\ requires representation of the same function in a
reference frame with a shifted origin compared to the original. This can be
handled using reexpansions of the basis functions, which for translations
along the $z$ axis gives%
\begin{equation}
S_{n}^{m}\left( \mathbf{r+i}_{z}\tau \right) =\sum_{n=0}^{p-1}\left(
S|R\right) _{n^{\prime }n}^{m}\left( \tau \right) R_{n^{\prime }}^{m}\left( 
\mathbf{r}\right) ,\quad R_{n}^{m}\left( \mathbf{r+i}_{z}\tau \right)
=\sum_{n=0}^{p-1}\left( R|R\right) _{n^{\prime }n}^{m}\left( \tau \right)
R_{n^{\prime }}^{m}\left( \mathbf{r}\right) ,  \label{shft0.1}
\end{equation}%
where $\left( S|R\right) _{n^{\prime }n}^{m}$ and $\left( R|R\right)
_{n^{\prime }n}^{m}$ are the coaxial reexpansion, or translation
coefficients (we do not show here the $\left( S|S\right) _{n^{\prime }n}^{m}$
$\left( \tau \right) $ coefficients which are equal to $\left( R|R\right)
_{n^{\prime }n}^{m}\left( \tau \right) $ coefficients. The problem then is
efficiently compute these coefficients.

The scheme for computation of the three-index coaxial translation
coefficients $\left( E|F\right) _{n^{\prime }n}^{m}$ ($E$ and $F$ can take
any value $S$ or $R$) proposed in \cite{Gumerov2004:SISC} and \cite%
{Gumerov2004:book} can be compacted to the following algorithm.

\begin{enumerate}
\item Compute zero-order sectorial coefficients for non-negative order $m$
(the coaxial coefficients depend only on $\left| m\right| $)%
\begin{eqnarray}
(S|R)_{n^{\prime }}^{0}(\tau ) &=&(-1)^{n^{\prime }}\sqrt{2n^{\prime }+1}%
h_{n^{\prime }}\left( k\tau \right) ,  \label{shft1} \\
(R|R)_{n^{\prime }}^{0}(\tau ) &=&(-1)^{n^{\prime }}\sqrt{2n^{\prime }+1}%
j_{n^{\prime }}\left( k\tau \right) ,  \notag \\
\left( E|F\right) _{n^{\prime }}^{m} &=&\left( E|F\right) _{n^{\prime
}m}^{m},\quad n^{\prime }=m,m+1,...,\quad E,F=S,R.  \notag
\end{eqnarray}

\item Compute all other order sectorial coefficients%
\begin{eqnarray}
b_{m+1}^{-m-1}\left( E|F\right) _{n^{\prime }}^{m+1} &=&b_{n^{\prime
}}^{-m-1}\left( E|F\right) _{n^{\prime }-1}^{m}-b_{n^{\prime }+1}^{m}\left(
E|F\right) _{n^{\prime }+1}^{m},  \label{shft2} \\
m &=&0,1,...,\quad n^{\prime }=m+1,m+2,...  \notag
\end{eqnarray}

\item Compute the tesseral coefficients using relation%
\begin{equation}
a_{n-1}^{m}\left( E|F\right) _{n^{\prime },n-1}^{m}-a_{n}^{m}\left(
E|F\right) _{n^{\prime },n+1}^{m}=a_{n^{\prime }}^{m}\left( E|F\right)
_{n^{\prime }+1,n}^{m}-a_{n^{\prime }-1}^{m}\left( E|F\right) _{n^{\prime
}-1,n}^{m}.  \label{shft3}
\end{equation}%
The latter recursion is applied to compute coefficients in the range%
\begin{equation}
m=0,...,p-1,\quad n=m,...,p-2,\quad n^{\prime }=n,...,2p-2-n,  \label{shft4}
\end{equation}%
using propagation in $n$ index%
\begin{equation}
\left( E|F\right) _{n^{\prime },n+1}^{m}=\frac{1}{a_{n}^{m}}\left[
a_{n-1}^{m}\left( E|F\right) _{n^{\prime },n-1}^{m}-a_{n^{\prime
}}^{m}\left( E|F\right) _{n^{\prime }+1,n}^{m}+a_{n^{\prime }-1}^{m}\left(
E|F\right) _{n^{\prime }-1,n}^{m}\right] .  \label{shft5}
\end{equation}

\item Compute coefficients for $n^{\prime }<n$ using symmetry%
\begin{equation*}
\left( E|F\right) _{n^{\prime }n}^{m}=(-1)^{n+n^{\prime }}\left( E|F\right)
_{nn^{\prime }}^{m}. 
\end{equation*}
\end{enumerate}

Numerical tests show that this scheme is instable at large $k\tau $. They
also show that steps 1) and 2) are stable, and the instability arises from
step 3). To stabilise the recursive computation we modify the algorithm as
follows. We use the existing scheme for $k\tau \lesssim 1$, and as a
relatively small number of translation coefficients need be computed, the
instability does not develop for reasonable computational accuracy (e.g.,
for double precision computing). For larger $k\tau \gtrsim 1$ instead of
propagation in $n$ index in step 3, see Eq. (\ref{shft5}), we use
propagation in the $n^{\prime }$ index,%
\begin{eqnarray}
\left( E|F\right) _{n^{\prime }+1,n}^{m} &=&\frac{1}{a_{n^{\prime }}^{m}}%
\left[ a_{n-1}^{m}\left( E|F\right) _{n^{\prime },n-1}^{m}-a_{n}^{m}\left(
E|F\right) _{n^{\prime },n+1}^{m}+a_{n^{\prime }-1}^{m}\left( E|F\right)
_{n^{\prime }-1,n}^{m}\right] ,\quad  \label{shift6} \\
m &=&1,...,p-1,\quad n^{\prime }=m+3,...,p-2,\quad n=m+1,...,n^{\prime
}-2.\quad  \notag
\end{eqnarray}%
For implementation some additional coefficients musr be computed, namely the
zonal and diagonal and subdiagonal coefficients. This introduces steps 3.1
and 3.2 which must be performed before step 3:

\begin{enumerate}
\item[3.1] Compute zonal coefficients using propagation in \ $n$-direction
for $\quad n=0,...,p-2,\quad n^{\prime }=n+1,...,2p-n-3$:\vspace*{-6pt} 
\begin{equation}
\left( E|F\right) _{n^{\prime },n+1}^{0}=\frac{1}{a_{n}^{0}}\left[
a_{n-1}^{0}\left( E|F\right) _{n^{\prime },n-1}^{0}-a_{n^{\prime
}}^{0}\left( E|F\right) _{n^{\prime }+1,n}^{0}+a_{n^{\prime }-1}^{0}\left(
E|F\right) _{n^{\prime }-1,n}^{0}\right] .  \label{shft7}
\end{equation}

\item[3.2] Compute diagonal and subdiagonal tesseral coefficients using
propagation in $m$ and in diagonals for $m=0,1,...,p-3,$ 
\begin{eqnarray}
\!\!\!\!\!\!(E|F)_{n+1,n}^{m+1} &=&\frac{1}{b_{n+1}^{-m-1}}\left(
b_{n+1}^{m}(E|F)_{n+1,n}^{m}+b_{n}^{-m-1}\left( E|F\right)
_{n,n-1}^{m}-b_{n}^{m}(E|F)_{n,n-1}^{m+1}\right) ,\,  \label{shft8} \\
n &=&m+2,...,p-2,  \notag \\
\!\!\!\!\!\!(E|F)_{n+1,n+1}^{m+1} &=&\frac{1}{b_{n+1}^{-m-1}}\left(
b_{n}^{m}(E|F)_{n+1,n-1}^{m+1}-b_{n+2}^{m}(E|F)_{n+2,n}^{m}+b_{n+1}^{-m-1}(E|F)_{n,n}^{m}\right) ,
\notag \\
\,n &=&m+1,...,p-2.  \notag \\
\!\!\!\!\!\!(E|F)_{n+2,n}^{m+1} &=&\frac{1}{b_{n+2}^{-m-1}}\left(
b_{n+1}^{m}(E|F)_{nn}^{m+1}-b_{n+1}^{m}(E|F)_{n+1,n+1}^{m}+b_{n+1}^{-m-1}(E|F)_{n+1,n-1}^{m}\right) ,
\notag \\
\,n &=&m+2,...,p-2.  \notag
\end{eqnarray}
\end{enumerate}

The reason for this modification of the original algorithm is in the
stability analysis based on the Courant-Friedrichs-Lewy (CFL) criterion \cite%
{Courant1967:IBM}. In step 3 of the original algorithm, the instability
manifests itself at large $n$ and $n^{\prime }$. We perform stability
analysis for the following recursion%
\begin{eqnarray}
\left( E|F\right) _{n^{\prime },n-1}^{m}-\left( E|F\right) _{n^{\prime
},n+1}^{m} &=&c\left[ \left( E|F\right) _{n^{\prime }+1,n}^{m}-\left(
E|F\right) _{n^{\prime }-1,n}^{m}\right] .  \label{shift9} \\
c &\sim &\frac{a_{n^{\prime }}^{m}}{a_{n}^{m}}\sim \frac{n}{n^{\prime }}%
\sqrt{\frac{n^{\prime 2}-m^{2}}{n^{2}-m^{2}}},\quad \left( n\gg 1,\quad
n^{\prime }\gg 1\right) ,  \notag
\end{eqnarray}%
Considering $n$ as time and $n^{\prime }$ as a spatial variable, then (\ref%
{shift9}) is like a finite-difference scheme for the wave equation with
celerity $-c$ and unit steps in the space and time. The CFL stability
condition is $\left| c\right| \leqslant 1,$ and so%
\begin{equation}
\frac{n}{n^{\prime }}\sqrt{\frac{n^{\prime 2}-m^{2}}{n^{2}-m^{2}}}\leqslant
1,\text{ or }n\geqslant n^{\prime }.  \label{shft10}
\end{equation}%
In the previous scheme the recursion was applied for computation of $\left(
E|F\right) _{n^{\prime },n+1}^{m}$ for $n^{\prime }\geqslant n$, which means
instability. The change in the direction of propagation fixes the problem.
The computation of the zonal and diagonal/subdiagonal coefficients (step 3.1
and 3.2 of the modified algorithm) are neutrally or asymptotically stable.
We implemented and tested the algorithm and found that it is stable and can
be used for computations of the reexpansion coefficients of large degree. We
also found that at low $k\tau $ the new algorithm exhibits a error due to
summation of terms of essentially different magnitude, leading to loss of
the precision. We thus propose the use of the original algorithm at low $%
k\tau $ and the modified algorithm for higher values of $k\tau $.

\section{Room Acoustic Simulations}
\label{Results} The BEM based on the theory presented above was implemented
using Open MP parallelization and validated. The results discussed below are
obtained on a workstation equipped with two Intel E5-2683 16 core processors
@2.1 GHz\ and 128 GB RAM.
\begin{figure}[tbh]
	\begin{center}
		\includegraphics[width=0.96\textwidth, trim=0.5in .5in 0.25in
		0.24in]{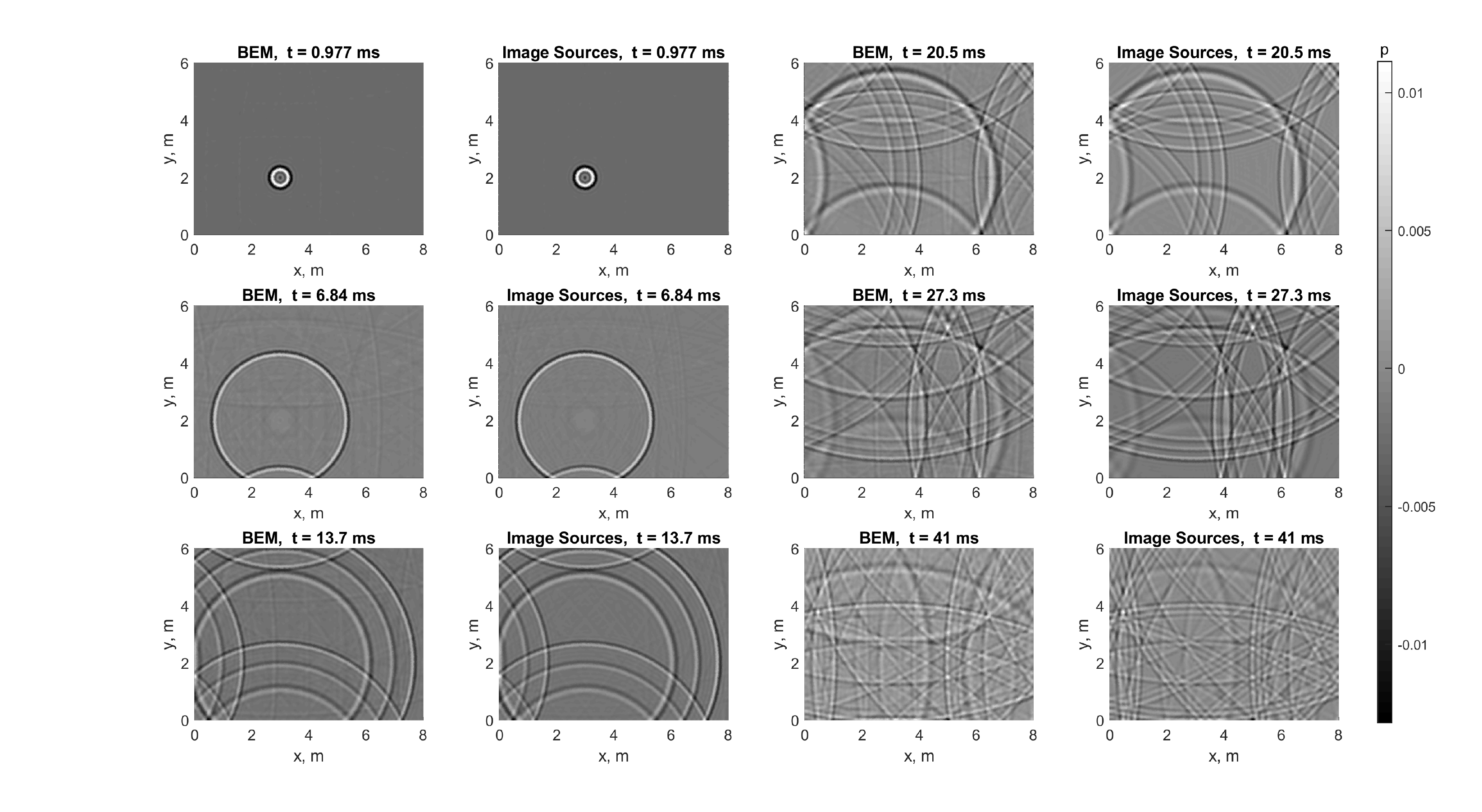}
	\end{center}
	\caption{A side-by-side comparison of the solutions at $z=1.75$ m obtained
		at different times using the direct BEM and the method of image sources for
		a pulse (\ref{an3}) from a source in a rectangular $8\times 6\times 3$ m
		room. The coordinates of the source are $\left( 3,2,1.75\right) $ m$.$ The
		reflection and absorption coefficients are $R=0.8$ and $\protect\alpha =0.36$
		, respectively. Computations performed for 128 frequencies and converted to
		the time domain using the FFT. }
	\label{Fig1}
\end{figure}

\begin{figure}[tbh]
	\begin{center}
		\includegraphics[width=0.96\textwidth, trim=0.5in .5in 0.25in
		0.24in]{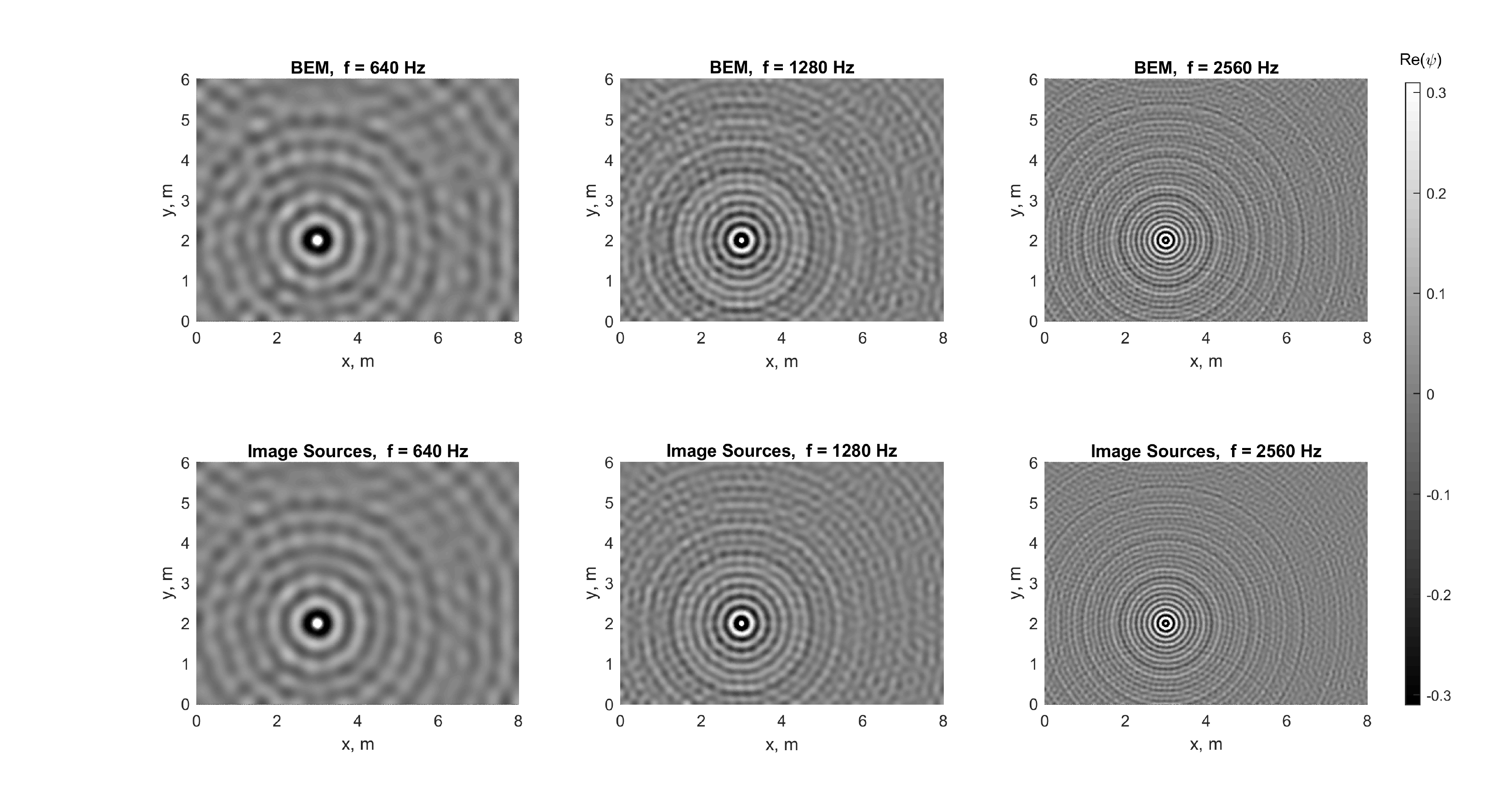}
	\end{center}
	\caption{ A comparison of the solutions described in Fig. \ref{Fig1} in the
		frequency domain at three different frequencies.}
	\label{Fig2}
\end{figure}
\begin{figure}[tbh]
	\begin{center}
		\includegraphics[width=0.96\textwidth, trim=0.5in .5in 0.25in
		0.24in]{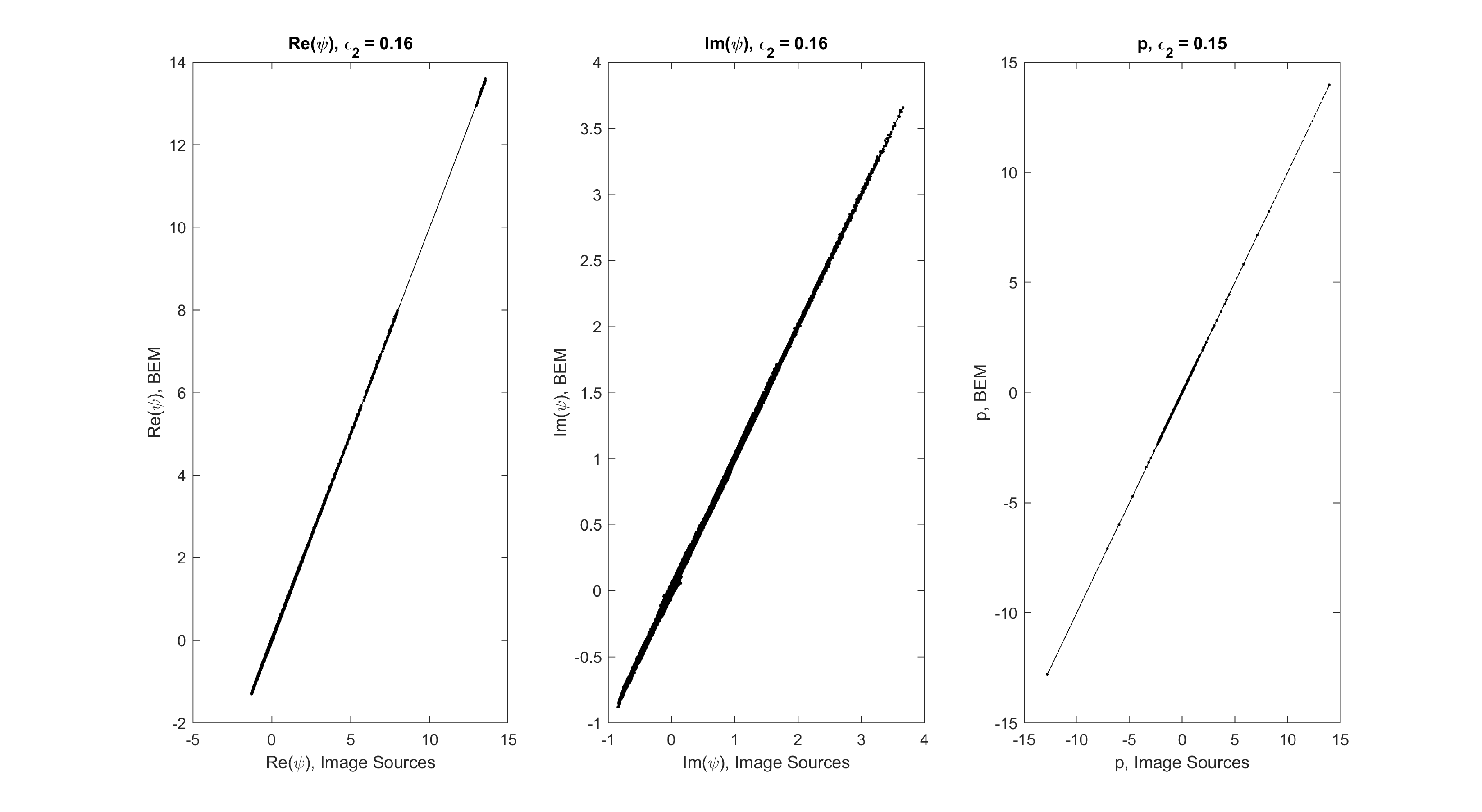}
	\end{center}
	\caption{Comparisons of the real and imaginary parts of the complex
		amplitude of the acoustic pressure obtained by the direct BEM and the image
		source method for all 128 frequencies and all evaluation points used for
		computing the images (on the left and in the center). A comparison for
		acoustic pressure at 256 different time instances and all evaluation points
		is shown on the rightmost graph. }
	\label{Fig3}
\end{figure}

\subsection{Comparisons with analytical solution}

For validation of the BEM we use analytical solution for a monopole source
located at $\mathbf{r=r}_{s}$ inside a rectangular room \cite{Allen1979:JASA}%
,%
\begin{equation}
p\left( \mathbf{r},t\right) =\frac{1}{\left| \mathbf{r-r}_{s}\right| }%
f\left( t-\frac{1}{c}\left| \mathbf{r-r}_{s}\right| \right) +\sum_{j}\frac{%
A_{j}}{\left| \mathbf{r-r}_{j}\right| }f\left( t-\frac{1}{c}\left| \mathbf{%
r-r}_{j}\right| \right) ,  \label{an1}
\end{equation}%
where $f\left( t\right) $ is a given function of time, and $\mathbf{r}_{j}$
are locations of the image sources of relative amplitude $A_{j}=R^{n_{j}}$,
where $R$ is the reflection coefficient and $n_{j}$ the index of the image
equal to the minimal number of reflections to obtain the image. For a
solution limited in time as $0\leqslant t\leqslant T_{\max }$, only sources
within radius $\left| \mathbf{r-r}_{j}\right| <cT_{\max }$ from evaluation
point $\mathbf{r}$ can be considered. Also the sources of intensity $A_{j}$
below some threshold, $A_{th}$, can be ignored (in our computations we used $%
A_{th}=10^{-7}$). \ For comparisons of solutions in the frequency domain we
also computed%
\begin{equation}
\psi \left( \mathbf{r},k\right) =G_{k}\left( \mathbf{r},\mathbf{r}%
_{s}\right) +\sum_{j}A_{j}G_{k}\left( \mathbf{r},\mathbf{r}_{j}\right) ,
\label{an2}
\end{equation}%
for each frequency. For comparisons with the analytical solution the latter
is generated for a room of dimensions $\left( l_{x},l_{y},l_{z}\right)
=\left( 8,6,3\right) $ m with wall reflection coefficient $R=0.8$. A
monopole sound source is placed at $\mathbf{r}_{s}=$ $\left( 3,2,1.75\right) 
$ m. For the function of time we used

\begin{equation}
f\left( t\right) =\delta _{f}\left( t-\Delta t\right) -\delta _{f}\left(
t+\Delta t\right) ,\quad \delta _{f}\left( t\right) =\left\{ 
\begin{array}{c}
1,\quad t=0 \\ 
0,\quad t\neq 0%
\end{array}%
\right.  \label{an3}
\end{equation}%
where $\delta _{f}$ is a finite version of the Dirac delta-function. In the
example shown below we used $T_{\max }=50$ ms, and $\Delta t=T_{\max }/256$.
All time values were interpolated to a time grid with step $\Delta t$ using
the closest neighbor value for shifted $f\left( t\right) $. We also computed
solution (\ref{an2}) for $128$ equispaced frequencies from 20 to 2560 Hz.
The results then were convolved with $f\left( t\right) $ using the FFT for
each evaluation point. We found good agreement between the frequency and
time domain computations, while some differences were explainable due to the
interpolation in the time domain. For further comparisons with the BEM we
used frequency domain computations.

The direct BEM for the Helmholtz equation with boundary conditions (\ref{sp7}%
), where $\zeta =9$ ($\alpha =1-R^{2}$ $=0.36$, see Eq. (\ref{sp4})) was
applied to get solution for all 128 frequencies. The room surface was
discretized with a mesh containing $N=1,806,532$ triangles with the maximum
edge length $2$ cm, which provided resolution at least 6.7 elements per
wavelength at the highest computed frequency. The field was evaluated at
241,542 grid points in a plain located at the height of the source ($z=1.75$
m) and covering the entire area of the room. The results of computations
were converted to the time domain as described above and compared with the
method of image sources.

Figure 1 shows comparisons of the results of the BEM\ and the image source
method in time domain. The time frames show excellent agreement and clear
physical interpretation of the results. At $t=0.977$ ms one can see a simple
spherical wavefront, which at the next time frame shown interacts with the
wall and the reflected wave is seen. Also, at $t=6.84$ ms one can see a dark
spot near the source location, which is due to the wave reflected from the
ceiling. The next time frame illustrates shows three wave fronts (triplet)
interacting with the walls. The wave front with the largest radius here is
due to the real source, while the two following wavefronts of smaller radii
are due to reflections from the ceiling and the floor, respectively. At $%
t=20.5$ ms one can see more reflections of the triplet of the wavefronts and
also a stronger single spherical front propagating from the source location
and reflected by the walls. This wavefront is a superposition of the waves
due to the secondary images located under the floor and above the sealing,
which arrival time to the imaging plane is the same. The last two figures
show more reflections. One can also observe the decay in the intensity of
multiple reflected waves.

Figure 2 illustrates the frequency domain computations at three selected
frequencies, which were used for synthesis of the time domain solution.
Again, one can see a good agreement between the computations using different
methods, while the wave pattern is quite complex due to multiple reflections.

Figure 3 compares point-by-point the analytical and BEM solutions in the
frequency domain for all 128 frequencies (the total number of data points
30,917,376 complex numbers) and in the time domain (61,834,752 real numbers)
for all 256 times. It is seen that the agreement is good, as ideally, all
points should fall on the line $y=x$. The relative $L_{2}$-norm error is
16\% for frequency domain and 15\% for time domain data.

We relate this difference to the boundary conditions used in the BEM and in
the method of images. Indeed, in the BEM we use impedance conditions (\ref%
{sp3}), which are conditions for the locally reacting surface, while the
method of images for absorbing walls, in fact, corresponds to the conditions
depending on the angle of incidence \cite{Kuttruff2017:book}. For specific
impedance $\zeta =0.9$ the deviation of the both types of conditions from
sound-hard boundary conditions can be of the order of 10\% or so. Note that
for time domain the errors are of the order of a few percents for relatively
small times when there are no reflections. Hence, 15\% difference between
the solutions should not be surprising, while the qualitative agreement of
the solution features is encouraging.

\subsection{More complex rooms}

\begin{figure}[tbh]
	\begin{center}
		\includegraphics[width=0.96\textwidth, trim=0.75in 1.5in 0.5in
		.4in]{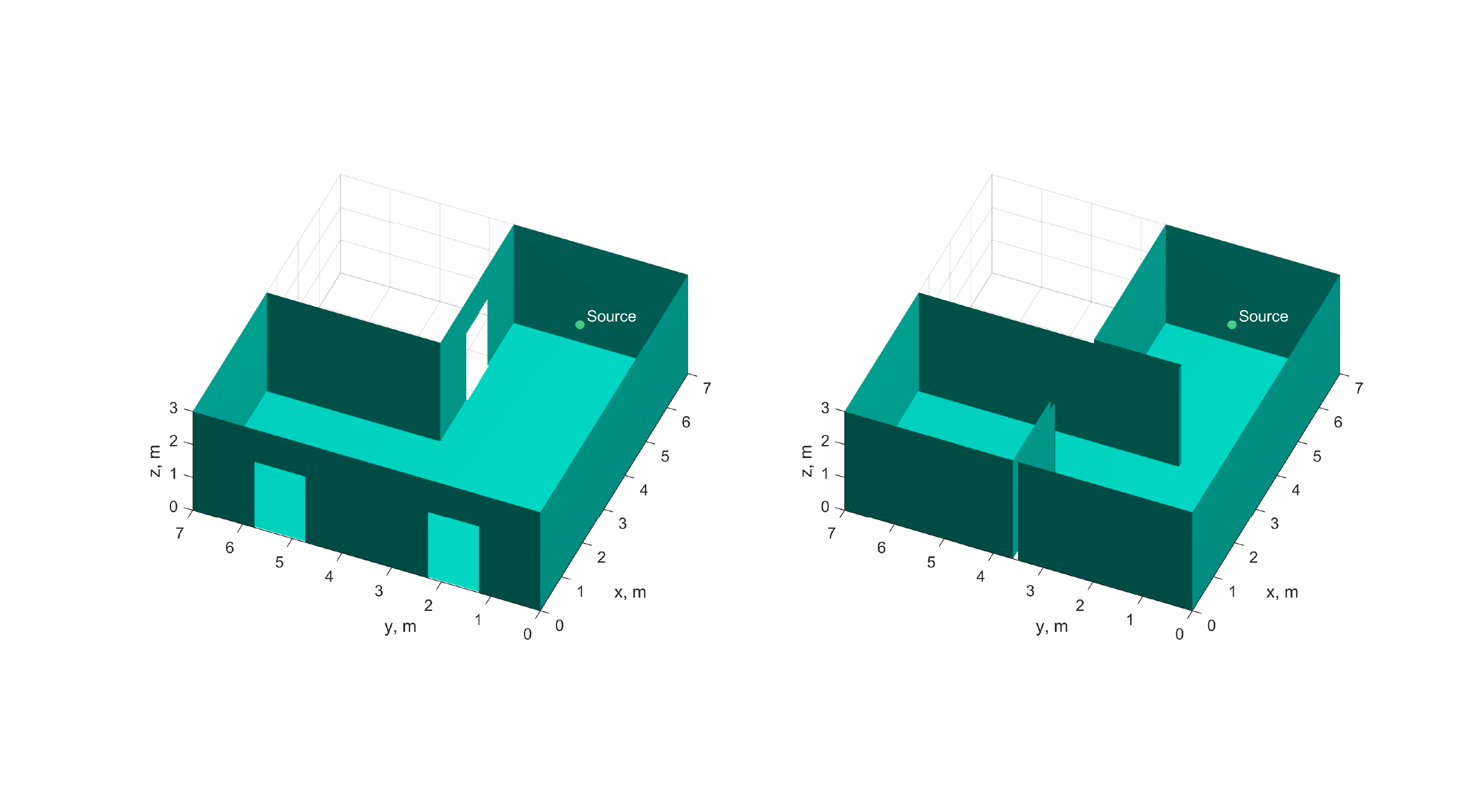}
	\end{center}
	\caption{Configurations of the L-shaped room used in simulations with open
		doors and with baffles (configurations of the room with closed doors/no
		baffles and without baffles/open doors are not shown). The room can be
		thought as a box of dimensions $7\times 7\times 3$ m, from which a box of
		dimensions $3.5\times 3.5\times 3$ m is extracted. The coordinates of the
		source are $\left( 6,1.75,1.75\right) $ m.}
	\label{Fig4}
\end{figure}
\begin{figure}[tbh]
	\begin{center}
		\includegraphics[width=0.96\textwidth, trim=0.75in .5in 0.5in .4in]{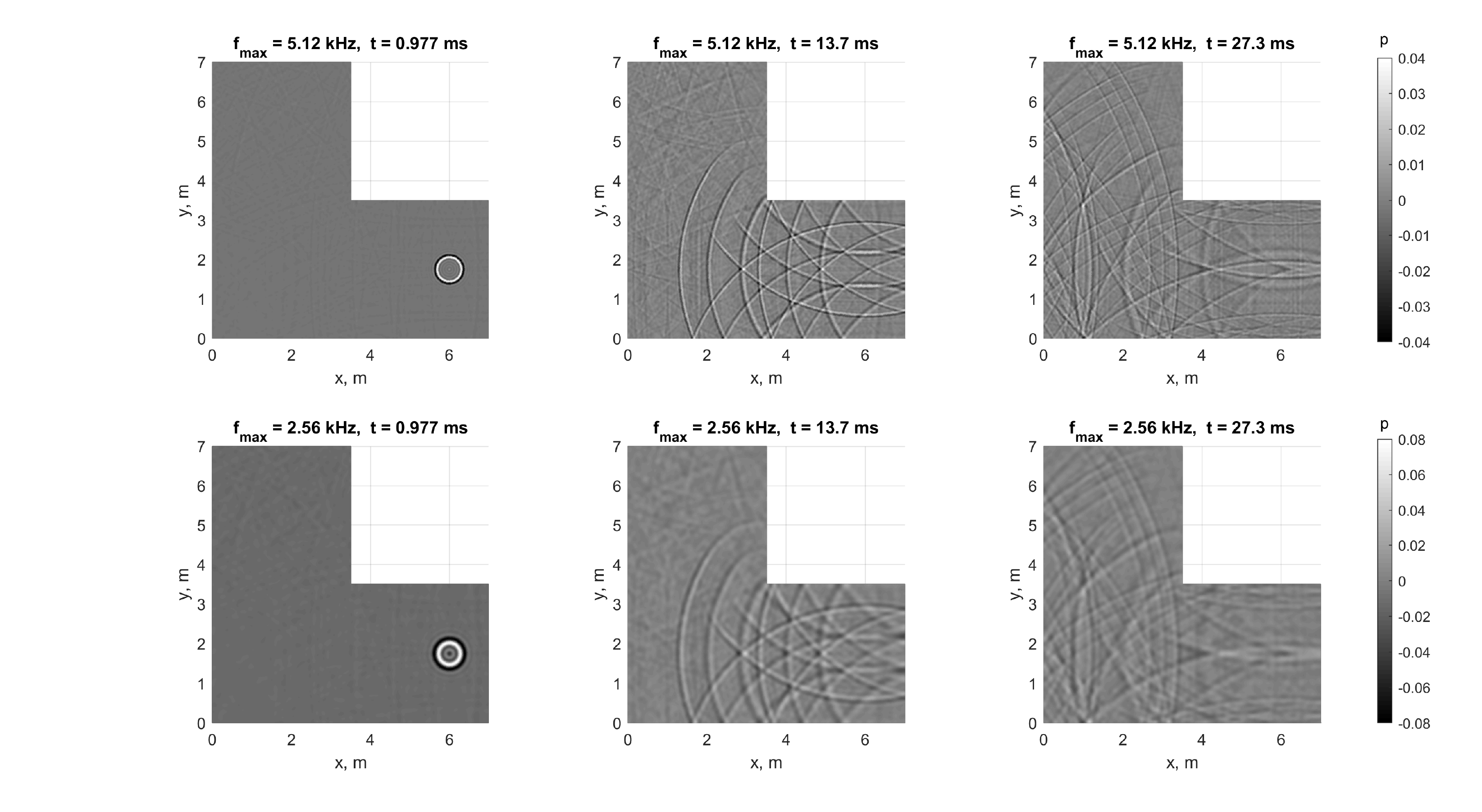}
	\end{center}
	\caption{Comparisons of the direct BEM computations at $z=1.75$ m and
		different times for the L-shaped room with closed doors using different
		spatial and temporal resolutions: a) mesh with 1 cm max element size and 256
		frequencies from 20 Hz to 5.12 kHz (the upper row); b) mesh with 2 cm max
		element size and 128 frequencies from 20 Hz to 2.56 kHz (the bottom row).
		The reflection and absorption coefficients are $R=0.9$ and $\protect\alpha %
		=0.19$. Patterns in time domain obtained using convolution of the pulse (\ref%
		{an3}) and frequency domain computations.}
	\label{Fig5}
\end{figure}

\begin{figure}[tbh]
	\begin{center}
		\includegraphics[width=0.96\textwidth, trim=0.75in 1.5in 0.5in
		.4in]{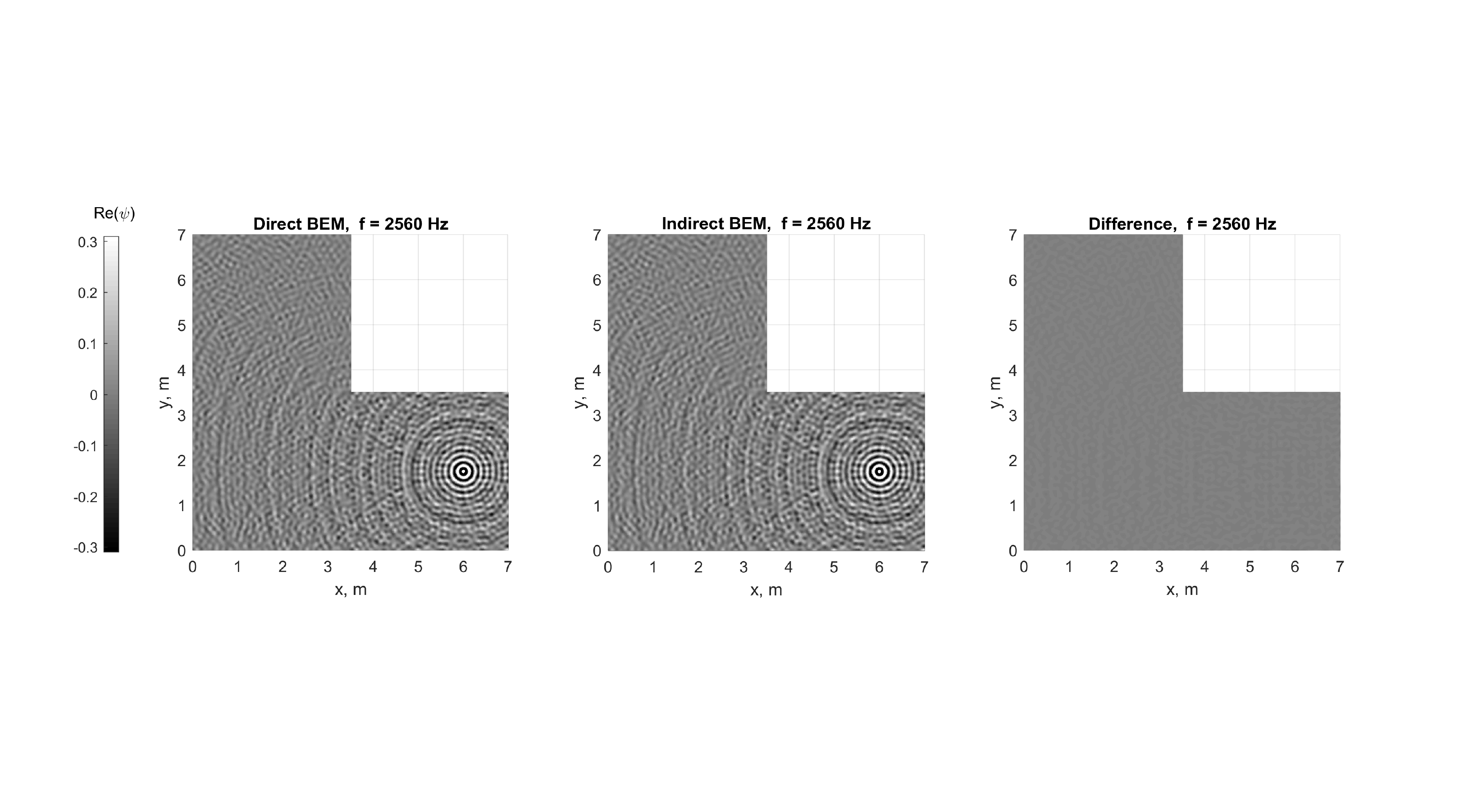}
	\end{center}
	\caption{Comparisons of the direct and indirect BEM computations at $z=1.75$
		m for the L-shaped room with closed doors at $f=2.56$ kHz, $R=0.9$ and $%
		\protect\alpha =0.19$.}
	\label{Fig6}
\end{figure}

\paragraph{L-shaped room}

Further, we conducted computations for more complex rooms, such as shown in
Fig. 4. Computations for the case shown in Fig. 4 on the left with the doors
closed are shown in Fig. 5. Two sets of computations for two different
meshes were performed. First, for the mesh with the maximum edge 1 cm
(6,310,256 triangles), and, second, for the mesh with the maximum edge 2 cm
(1,579,544 triangles). The first mesh was used for computations for
frequencies up to 5.12 kHz ($kD$ $\approx $ $1137$), while the second mesh
was used for computations for frequencies up to 2.56 kHz ($kD$ $\approx $ $%
569$). This corresponds to the mesh resolution of 6.7 elements per
wavelength for the highest frequency. In both cases the reflection
coefficient was set to a constant $R=0.9$ ($\alpha =0.19,$ $\zeta =19$). The
time domain data are obtained as in the case for the rectangular room using
the FFT for 256 and 128 frequencies equispaced from 20 Hz to the maximum
frequency. Note that the case for a finer mesh is evaluated on a finer grid
consistent with the mesh vertices.

Figure 5 compares time domain computations for two cases for pulses
described by Eq. (\ref{an3}). The comparison shows agreement between the two
data sets, while computations on finer mesh provide sharper pictures. We
also compared computations using the direct and indirect methods. Figure 6
compares side-by-side simulations using these methods in the frequency
domain for 2.56 kHz. It is hard to find any difference visually, while the
difference of about 4.5\% in the relative $L_{2}$-norm exists. Such level of
the errors can be related to the BEM itself (discretization, approximation
of far-field integrals, etc.). The comparison can be used as
cross-validation of the direct and indirect methods.

\paragraph{Room with openings}

Opening in a room such as open doors, windows or ventilation outlets can be
modeled in two ways. The more natural way is to use the indirect BEM, which
can handle open surfaces. Another way is to use the direct BEM with proper
boundary conditions imposed on the regions of openings. For example, one can
assign to the impedance boundary condition (\ref{sp3}) with specific
impedance $\zeta =1$, which corresponds to the reflection coefficient $R=0$
(absorption $\alpha =1$, see Eq. (\ref{sp4})).
\begin{figure}[tbh]
	\begin{center}
		\includegraphics[width=0.96\textwidth, trim=0.75in 0.5in 0.5in
		.4in]{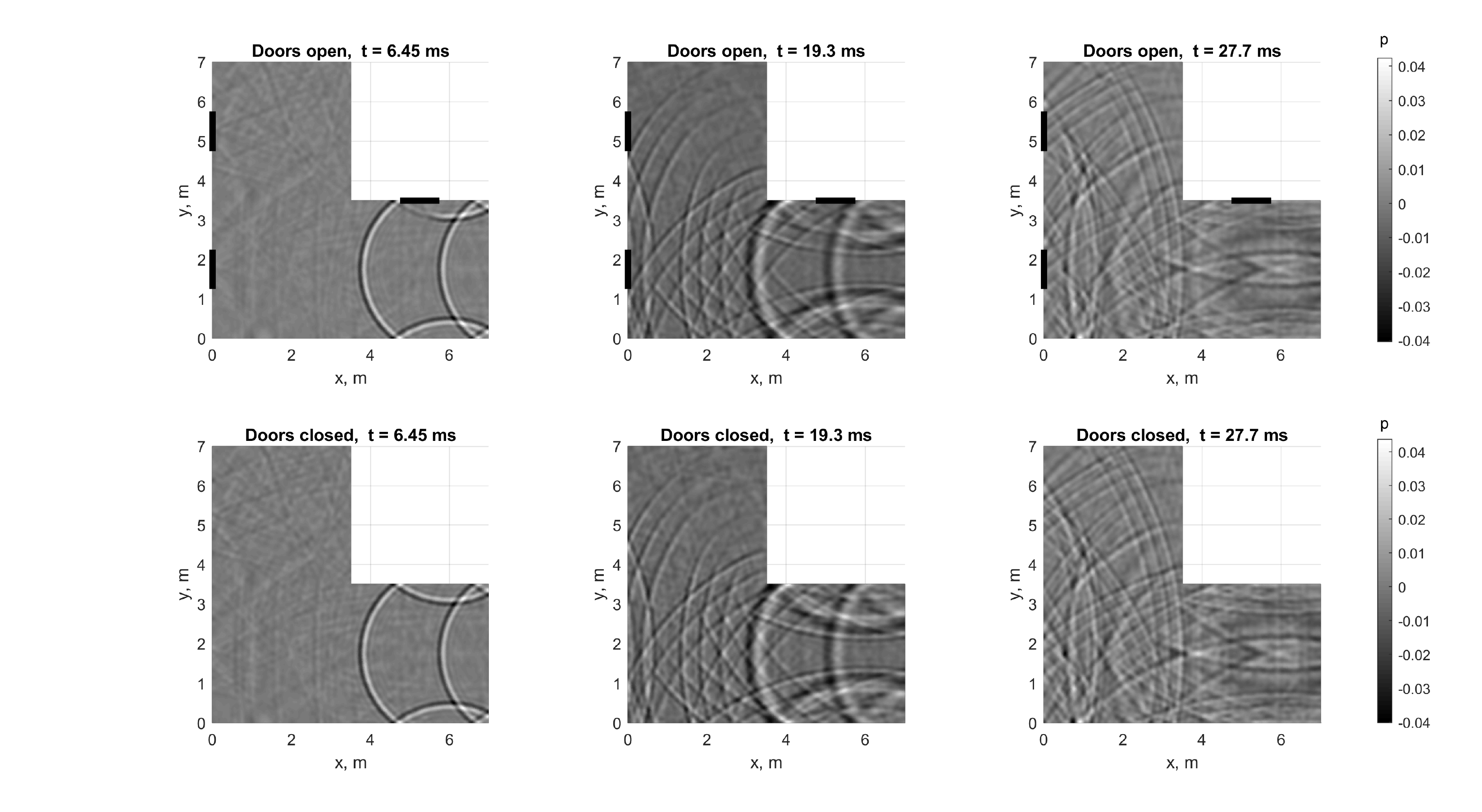}
	\end{center}
	\caption{Comparisons of the direct BEM computations at $z=1.75$ m and
		different times for the L-shaped room with closed and open doors. The
		reflection and absorption coefficients are $R=0.9$ and $\protect\alpha %
		=0.19. $ The doors are modelled as regions with $R=1$ and $\protect\alpha =0.
		$ }
	\label{Fig7}
\end{figure}

\begin{figure}[tbh]
	\begin{center}
		\includegraphics[width=0.96\textwidth, trim=0.75in 1.5in 0.5in
		.4in]{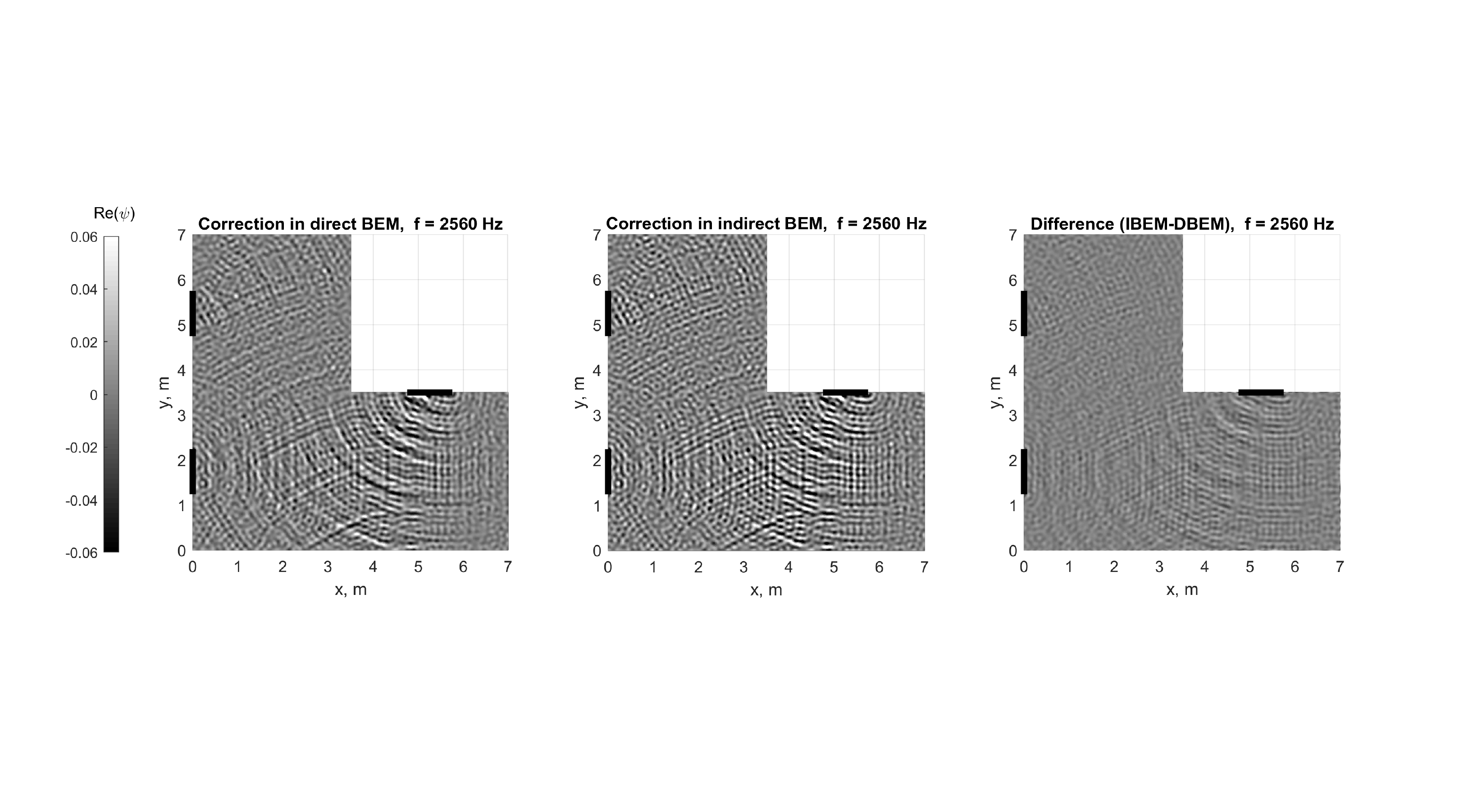}
	\end{center}
	\caption{Correction to the closed door L-shaped room simulations to account
		for open doors when using the direct and indirect BEM at $z=1.75$ m, $f=2.56$
		kHz, $R=0.9$ and $\protect\alpha =0.19.$ In the direct BEM the doors are
		modelled as regions with $R=1$ and $\protect\alpha =0.$ The chart on the
		right shows the difference between the indirect and direct BEM computations. 
	}
	\label{Fig8}
\end{figure}

\begin{figure}[tbh]
	\begin{center}
		\includegraphics[width=0.96\textwidth, trim=0.75in 0.5in 0.5in
		.4in]{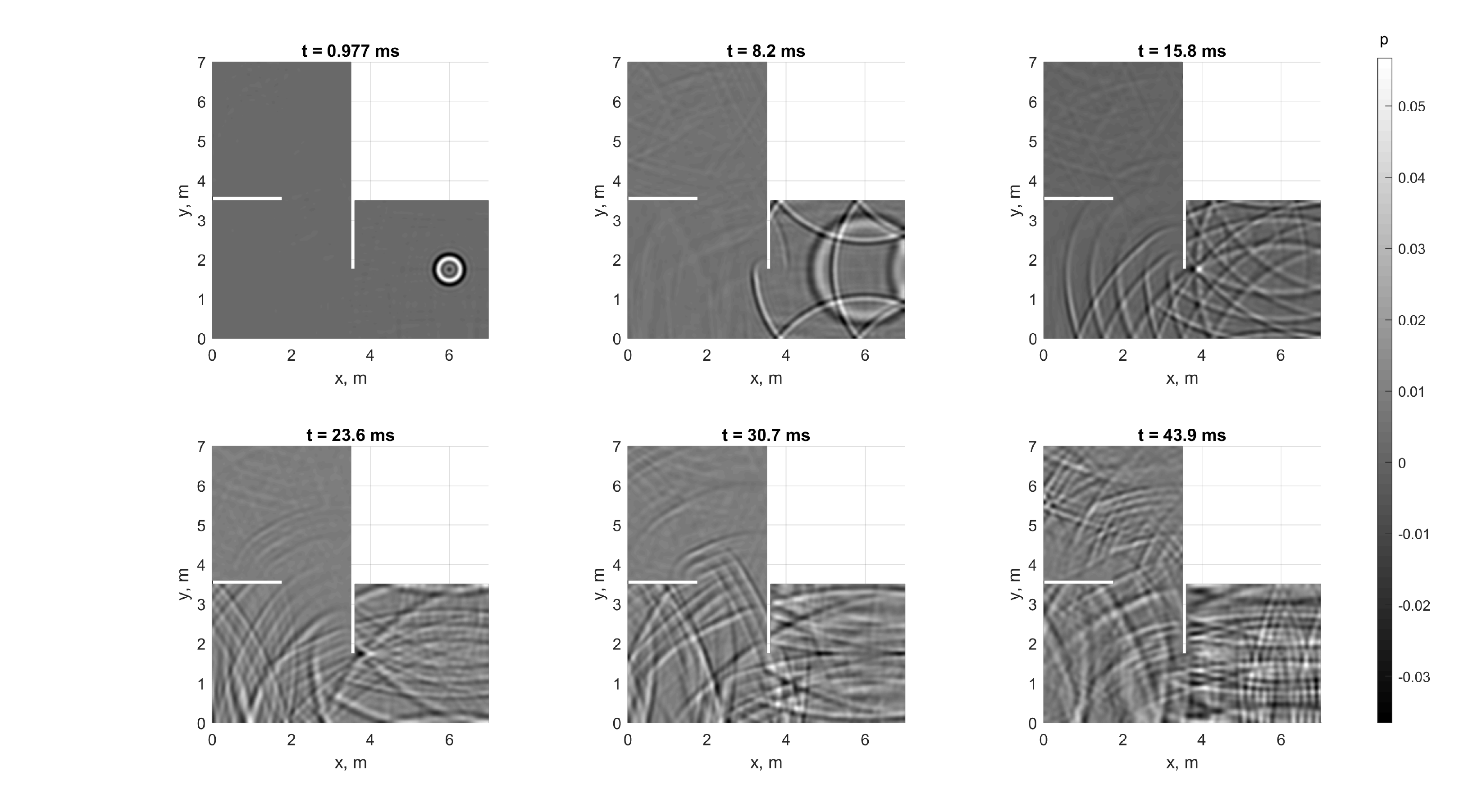}
	\end{center}
	\caption{Direct BEM computations of pulse (\ref{an3}) propagation in the
		L-shaped room with baffles (closed doors) at $z=1.75$ m and different times.
		The reflection and absorption coefficients are $R=0.9$ and $\protect\alpha %
		=0.19$ everywhere but on th baffles at which $R=0$ and $\protect\alpha =1$.
		Computations performed for 128 frequencies from 20 Hz to 2.56 kHz and
		converted to the time domain using the FFT.}
	\label{Fig9}
\end{figure}

Figure 7 illustrates computations for the same 2 cm resolution (maximum
frequency 2.56 kHz ) as in the case illustrated in Fig. 5 ($R=0.9$ on the
entire surface). Here the cases with the doors closed and open are compared
and overall difference of 20\% in the relative $L_{2}$-norm found (all
receiver points over all times). While the wave patterns look similar for
these two cases, one can notice some differences near the door regions, as
the reflected waves in the case of closed doors are stronger. This
difference is better seen as we plot the difference between the cases of
open an closed doors shown in Fig. 8 in the frequency domain. Here one can
see that the difference is equivalent to the effect of sound sources
distributed over the opening areas. Also, on the same figure we displayed
the difference between the cases of open and closed doors computed using the
indirect BEM. It is seen that the wave patterns for the direct and indirect
BEM simulations are the same. Nonetheless, there exists a difference of
about 8\% between computations obtained by the two methods, which is larger
than the expected BEM errors. This difference is also illustrated in Fig. 8
and is associated with the fact that boundary conditions (\ref{sp3}) with
specific impedance $\zeta =1$ are not fully absorbing (they are fully
absorbing only for the waves of normal incidence). So, some discrepancy
between the direct BEM with conditions (\ref{sp3}) and indirect BEM (the
latter corresponds to truly full absorption) should exist and we observe it
in our computations. Note here, that for modeling of openings some other
boundary conditions and methods instead of Eq. (\ref{sp3}) can be used, for
example, one can consider utilization of perfectly matched layers, which is
common for computations of wave equations in open domains.

\paragraph{Room with baffles}

\begin{figure}[tbh]
	\begin{center}
		\includegraphics[width=0.96\textwidth, trim=0.75in 1.5in 0.5in
		.4in]{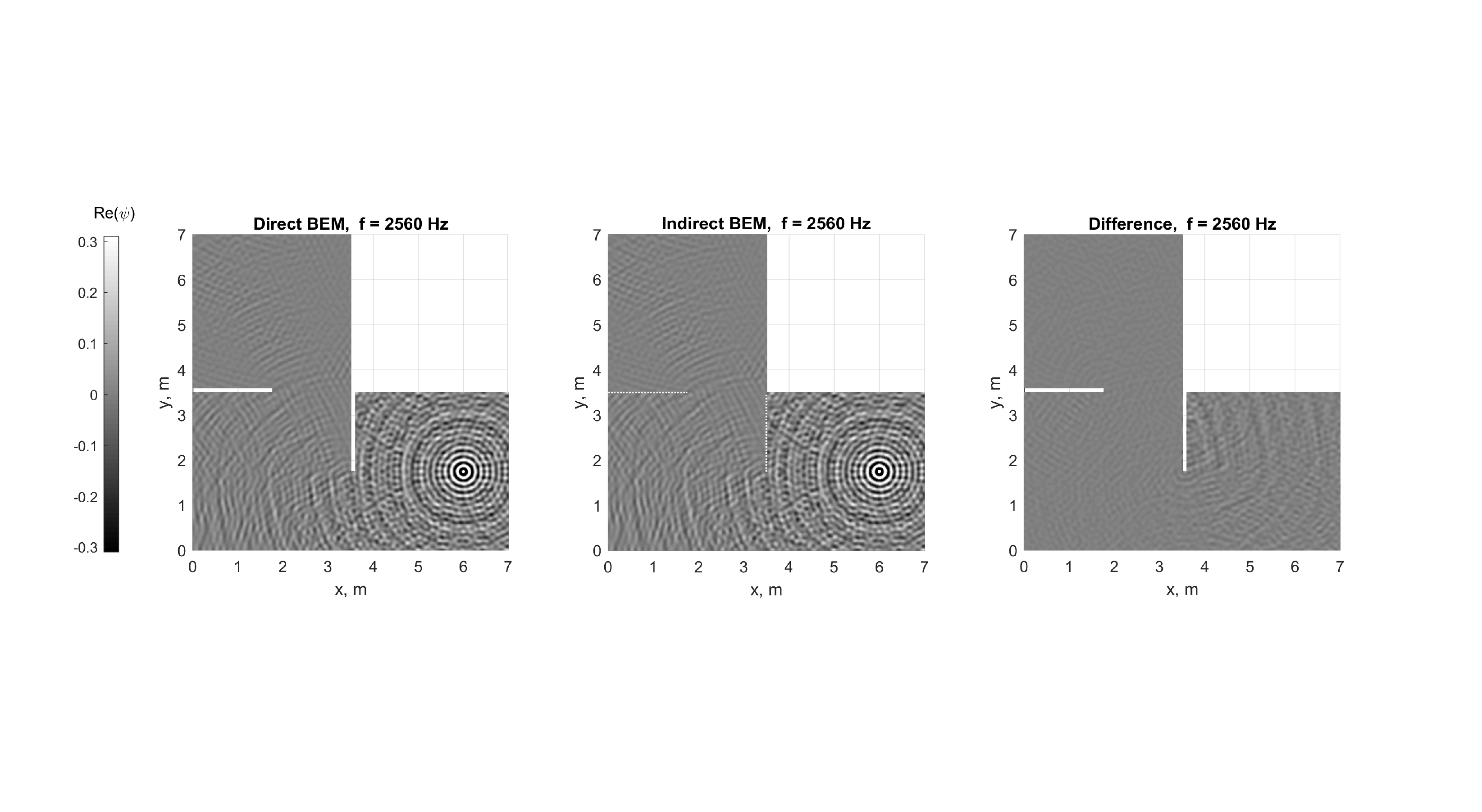}
	\end{center}
	\caption{Comparisons of the direct and indirect BEM computations at $z=1.75$
		m for the L-shaped room with baffles (closed doors) at $f=2.56$ kHz, $R=0.9$
		and $\protect\alpha =0.19$ ($R=0$ and $\protect\alpha =1$ on the baffle).
		For the direct BEM, the thickness of the baffle is 7 cm.}
	\label{Fig10}
\end{figure}

The presence of baffles drastically changes the acoustic field in the room.
For computational examples we used the L-shape room with closed doors shown
in Fig. 4. The reflection coefficient is $R=0.9$ for the entire room, while
we assumed that the baffles are fully absorbing, $R=0$. Again, two methods
were used for computations and analysis, the direct and indirect BEM. In the
direct BEM the thickness of baffles is 7 cm$,$ while in the indirect method
the baffles are treated as infinitely thin objects.

Figure 9 illustrates several frames in the time domain obtained using direct
BEM executed for a 2 cm resolution mesh and maximum frequency 2.56 kHz. The
frame at $t=8.2$ ms shows the principal wave and its reflections from the
walls (thin wave fronts) and from the ceiling (thick wave fronts). It also
shows interaction of the main wave with the baffle, as only part of the wave
goes through the opening, while rather weak reflection is observed from the
baffle itself. This also confirms that the baffle does not absorb the wave
perfectly as the angular dependence of the reflection coefficient is ignored
in the model. However, for the frames at later stages the reflections from
the baffles are not observed clearly, which also can be related to
substantial damping of such reflections from the image sources. Also, we
analyzed wave patterns for several frequencies and could not see visually
substantial reflections from the baffles. Examples of such patterns obtained
using the direct and indirect BEM are shown in Figure 10 for frequency 2.56
kHz. Note that visually computations for the direct and indirect methods are
very similar, but, in fact, the difference of 16\% in the relative $L_{2}$%
-norm exists. Such difference cannot be explained neither by the BEM
accuracy, nor the models of reflection, as for both methods we used the same
impedance boundary conditions with $\zeta =1$ on the baffles.

The explanation comes from the consideration of the difference of the fields
computed by the two methods also illustrated in Fig. 10. One can clearly see
some pattern reminding a reflection from the baffle closest to the source.
We relate this to the thickness of the baffle. Indeed, at 2.56 kHz the
wavelength is 13.4 cm, so the baffle is approximately half-wavelength thick.
By satisfying boundary conditions on the baffles we should observe that the
crests and the troughs of the reflected waves near the offset baffle points
should be exchanged and amplified in the difference figure. Note that we
could not see any regular patterns in Fig. 6 illustrating the difference
between the two methods for a room without baffles.

\subsection{Performance}
\begin{figure}[tbh]
	\begin{center}
		\includegraphics[width=0.96\textwidth, trim=0.75in 0.5in 0.5in
		.4in]{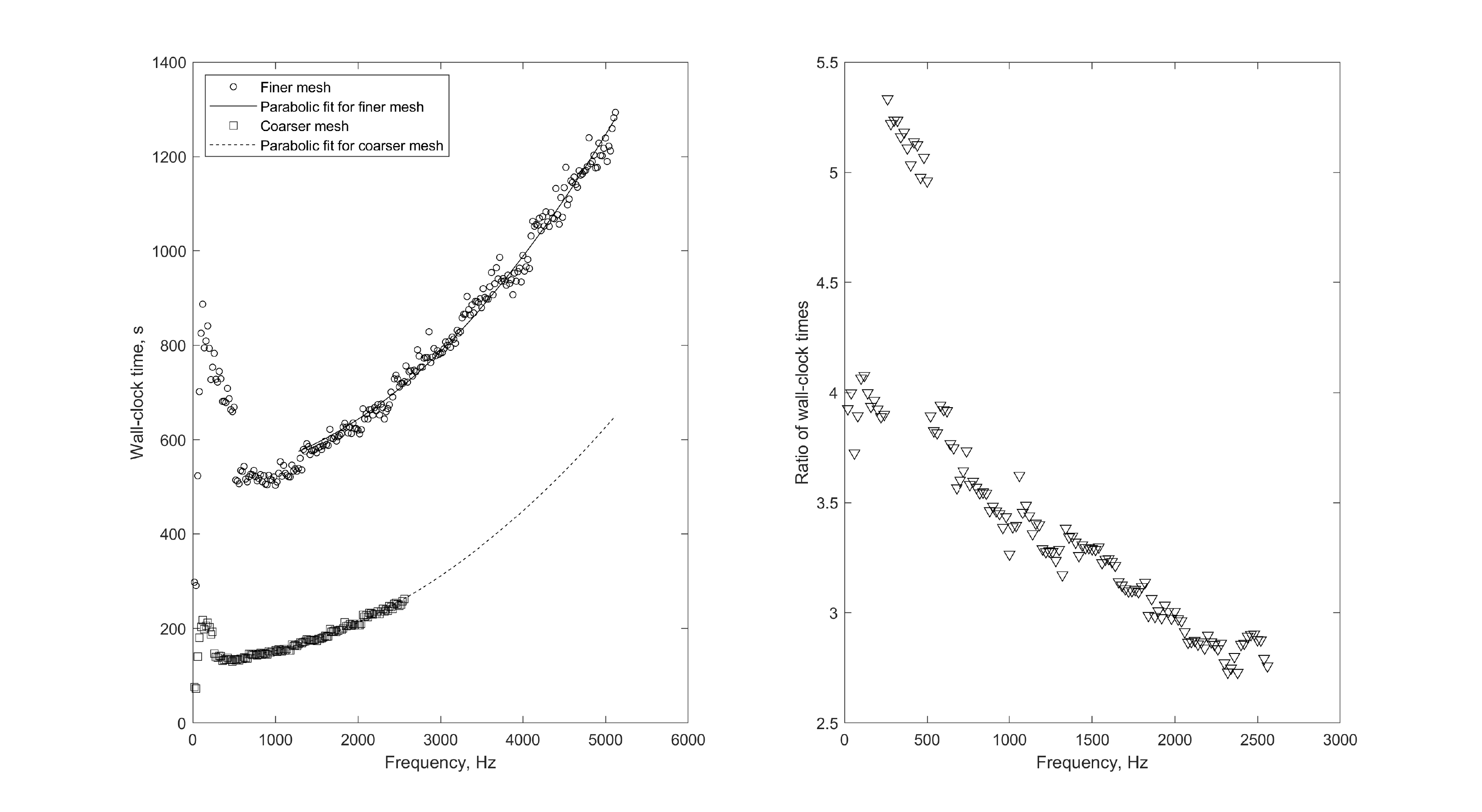}
	\end{center}
	\caption{On the left: the wall-clock times required for direct BEM
		simulations for the L-shaped room with closed doors at different frequencies
		and different mesh sizes, $N=6,310,256$ triangles (maximum edge 1 cm) and $%
		N=1,579,544$ triangles (maximum edge 2 cm), $R=0.9$ and $\protect\alpha %
		=0.19 $. The continuous lines show parabolic fits On the right: the ratio of
		the wall-clock times for the higher and lower mesh resolutions.}
	\label{Fig11}
\end{figure}

\begin{figure}[tbh]
	\begin{center}
		\includegraphics[width=0.96\textwidth, trim=0.75in .5in 0.5in
		.4in]{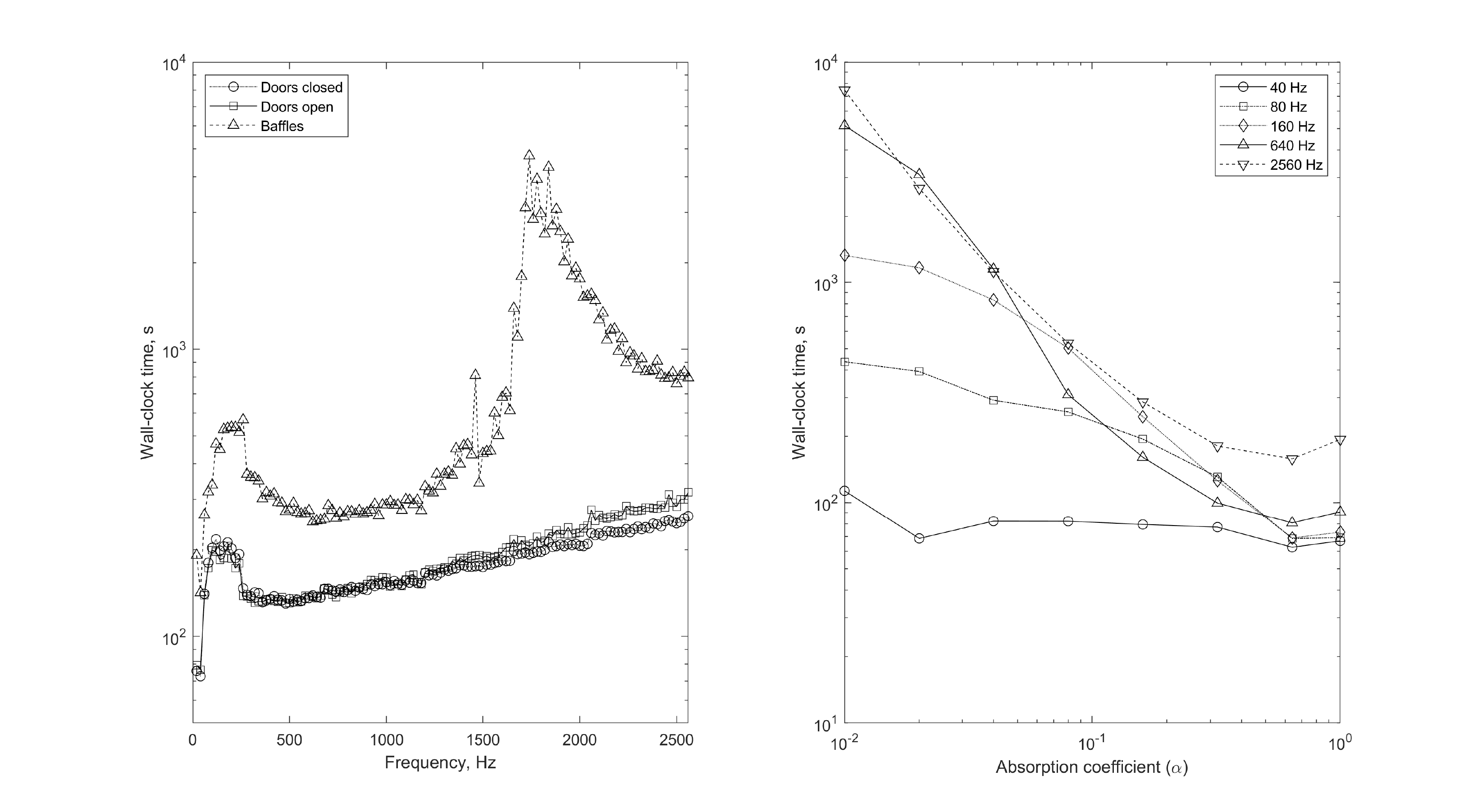}
	\end{center}
	\caption{The influence of different factors on the time required for
		solution. On the left: the wall-clock time vs frequency required for the
		direct BEM solution of different configurations of the L-shaped room. On the
		right: the dependence of the direct BEM solution time for the L-shaped room
		with closed doors on the absorption coefficient $\protect\alpha $.}
	\label{Fig12}
\end{figure}

Some remarks can be made concerning the performance of the solvers. In most
cases we computed using the direct BEM formulation (\ref{dbem4}) provided
much better convergence than (\ref{dbem5}). Figure 11 illustrates wall-clock
times using two different meshes for L-shape room (for the case illustrated
in Fig. 4). At fixed frequency the FMM is scaled almost linearly with
respect to $N$ (more precisely, as $O\left( N\log ^{m}N\right) $, $m>0$). As
the mesh with resolution 1 cm contains 4 times more elements than the mesh
with resolution 2 cm, ideally the solution time should be 4 times larger for
the finer mesh. The ratio of the observed wall-clock times required for
solution varies around this value (as the solution itself depends also on
the number of matrix-vector calls in the iterative processes). Figure 11
also shows that at higher frequencies the wall-clock times, $w$, for the
present algorithm are fitted well with a fit $w\sim a+bf^{2}$, where $f$ is
the frequency. Theoretically at fixed $N$ and constant number of iterations
the scaling for the present algorithm should be $w\sim a+bf^{2}+cf^{3}$ with
relatively small constant $c$ \cite{Gumerov2009:JASA}.

Figure 12 illustrates various factors affecting the performance. The figure
on the left panel shows that the performance of the direct BEM for the
L-shape room with closed and open doors approximately the same. However, the
times for computation of the room with the baffles can be several times
larger. Note that the computation times are substantially affected by
settings of the flexible GMRES preconditioned with the low accuracy FMM
GMRES, such as the tolerance for convergence, the dimensionality of the
Krylov subspaces, the number of iterations for the preconditioner, and the
accuracy of the FMM used in the main and in the preconditioning loops \cite%
{Gumerov2009:JASA}. The data plotted in the figure are obtained by taking
the fastest times from several runs with different parameter settings. While
this can be considered as some manual optimization, a systematic study of
optimal parameters for the iterative solvers is required. The figure on the
right panel of Fig. 12 shows the effect of the absorption coefficient $%
\alpha $ on the wall-clock time. It is seen that for low frequencies like 40
Hz the computation time just slightly depends on $\alpha $. Such dependence
is much stronger for larger frequencies and the wall-clock times at $\alpha
\sim 0.01$ can be an order of magnitude larger than those at $\alpha \sim 1$%
. This figure also shows that at fixed $\alpha $ the frequency dependence of
the wall-clock time can be non-monotonic. Note also that the wall-clock
times for indirect BEM at large enough $\alpha $ are of the same order of
magnitude as for the direct BEM, while for smaller $\alpha $ they can be
much larger. Moreover, one may need to use some additional or different
preconditioners in this case (we used a simple diagonal preconditioner (a
constant diagonal added to the preconditioner matrix)).

Finally, we note that the present approach has a great opportunity for easy
distributed parallelization as the solution at different frequencies can be
obtained independently on different computing nodes of a cluster. For
example, Fig. 11 shows that for any of 128 frequencies computations on the 2
cm resolution mesh took less than 300 seconds. This means that a cluster of
128 nodes similar to our workstation can solve the entire problem for 300
seconds or so (in fact, a smaller number of nodes can be used for this and
other cases as some nodes can perform two or more jobs; also, some
optimization problems for job distributions can be considered).

\section{Conclusion}

Our computations, comparisons, and analysis show that the boundary element
methods can be used for simulation of room acoustics for small rooms $V\sim
150$ m$^{3}$. There are a number of advantages for the use of the BEM for
room acoustics as this method is free from a number of deficiencies of the
methods based on geometric acoustics approximation, provides a comprehensive
information about the sound fields inside a room, and enables accounting for
frequency dependent spatially varying absorption. Comparisons with the image
source method for a rectangular room shows that the neglect of the angular
dependence of the impedance generates a relatively small error in the
solution of the order of 15\%. Computations for an L-shaped room show that
the effects of openings such as doors can be of the order of 20\%. Modeling
of baffles should be done with care as the thickness of baffles can cause
15\% difference at wavelengths comparable with the baffle thickness.

Our numerical experiments also show that a proper form of the boundary
integral equation should be selected for the direct method. The difference
in performance heavily depends on this formulation. As a rule, for low
absorptions the indirect BEM shows much slower convergence than the direct
BEM. Furthermore, there are many factors affecting the BEM performance, such
as the presence of baffles and the value of the absorptions coefficient. It
is also important to find some optimal settings for iterative solver. Tuning
of the solver parameters can improve performance several times. While in
some cases the convergence of the BEM is excellent and the problem for the
entire frequency range can be solved on computational clusters for 5 min or
so, there exist cases of slow convergence, which require more detailed
studies, and probably design of good preconditioners.

\section*{Acknowledgments}
	This research was conducted at VisiSonics Corporation.

\nocite{*}
\bibliographystyle{alpha}
\bibliography{RoomAcoustics.bib}
\end{document}